\documentclass[10pt]{amsart}
\usepackage[ansinew]{inputenc}
\usepackage[all]{xy}
\usepackage{amsmath,latexsym,amssymb,verbatim,amsthm}

\input arrow.tex

\newtheorem{lema}{Lemma}[section]
\newtheorem{teo}[lema]{Theorem}
\newtheorem{pro}[lema]{Proposition}

\newtheorem{cor}[lema]{Corollary}

\newtheorem{teo*}{Theorem}

\begin{document}

\title{Characterization of quasi-coherent modules that are module schemes}

\author{Amelia \'{A}lvarez}
\address{Departamento de Matem\'{a}ticas, Universidad de Extremadura,
Avenida de Elvas s/n, 06071 Badajoz, Spain}
\email{aalarma@unex.es}

\author{Carlos Sancho}
\address{Departamento de Matem\'{a}ticas, Universidad de Salamanca,
Plaza de la Merced 1-4, 37008 Salamanca, Spain}
\email{mplu@usal.es}

\author{Pedro Sancho}
\address{Departamento de Matem\'{a}ticas, Universidad de Extremadura,
Avenida de Elvas s/n, 06071 Badajoz, Spain} \email{sancho@unex.es}





\date{June, 2007}

\maketitle


\section*{Introduction}

Let $R$ be a commutative ring with unit. All functors we consider
are functors over the category of commutative $R$-algebras. Given
an $R$-module $E$, we denote by ${\bf E}$ the functor of
$R$-modules ${\bf E}(B) := E \otimes_R B$. We will say that ${\bf
E}$ is a quasi-coherent $R$-module. If $E$ is an $R$-module of
finite type we will say that ${\bf E}$ is a coherent $R$-module.
Given $F$, $H$ functors of $R$-modules, ${\bf Hom}_R (F,H)$
will denote the functor of $R$-modules $${\bf Hom}_R (F,H)(B)=
{\rm Hom}_B(F_{|B}, H_{|B})$$ where $F_{|B}$ is the functor $F$
restricted to the category of commutative $B$-algebras. It holds
that ${\rm Hom}_R ({\bf E}, {\bf E}') = {\rm Hom}_R (E,E')$. The
category of $R$-modules is equivalent to the category of
quasi-coherent $R$-modules (\cite[1.12]{Amel}).

We denote $F^* := {\bf Hom}_R(F, {\bf R})$. We will say that ${\bf
E}^*$ is an $R$-module scheme.

The $R$-module functors that are essential for the development of
the theory of the linear representations of an affine $R$-group
are the quasi-coherent $R$-modules and the $R$-module schemes
(\cite{Amel}). The aim of this paper is to study when a
quasi-coherent $R$-module is an $R$-module scheme. We will prove
that it is equivalent to giving a characterization of projective
$R$-modules of finite type.

The main result we are going to use is the following proposition.

\begin{pro}\label{0.1} \cite[1.8]{Amel}
Let $E$, $E'$ be $R$-modules. Then:
$${\bf Hom}_R({\bf E}^*, {\bf E}') = {\bf E} \otimes_R {\bf E}'.$$
%
Two immediate consequences are:
\begin{enumerate}
\item[(a)] ${\bf E}^{**} = {\bf E}\; ($\cite[1.10]{Amel}$)$.

\item[(b)] If $E$ is a projective module, the image of a morphism
${\bf E}^* \to {\bf V}$ is a coherent module
$($\cite[4.5]{Amel}$)$.
\end{enumerate}
\end{pro}

\section{Characterization of quasi-coherent modules that are module schemes}

\begin{teo}
Let $E$ be an $R$-module. Then, ${\bf E}$ is an $R$-module scheme if and only if $E$ is a projective module of finite type.
\end{teo}

\begin{proof}
$\Rightarrow$) Assume ${\bf E} = {\bf V}^*$. For every $R$-module $E'$, $${\rm Hom}_R (E, E')= {\rm Hom}_R ({\bf E}, {\bf E}') = {\rm Hom}_R ({\bf V}^*, {\bf E}') \stackrel{\text{\ref{0.1}}}{=} V \otimes_R E' .$$ $E$ is a projective module because ${\rm Hom}_R(E, -) = V \otimes_R -$,
which is exact on the right.

Since ${\bf V} \stackrel{\text{\ref{0.1}(a)}}{=} {\bf V}^{**} =
{\bf E}^*$, we have that $V$ is a projective module. The image of
the isomorphism ${\bf V}^* = {\bf E}$ is coherent (\ref{0.1} (b)).
Then, $E$ is an $R$-module of finite type.

$\Leftarrow$) Let us consider an epimorphism from a finitely
generated free module $L$ to $E$, $L \to E$. Taking dual we have
an injective morphism ${\bf E}^* \hookrightarrow {\bf L}^*$. ${\bf
L}^*$ is isomorphic to a coherent module. The image of the
morphism ${\bf E}^* \hookrightarrow {\bf L}^*$, which is ${\bf
E}^*$, is coherent (\ref{0.1}(b)). Then, ${\bf E} = {\bf E}^{**}$ is
a module scheme.
\end{proof}

\begin{cor}
An $R$-module $E$ is projective of finite type if and only if
$${\rm Hom}_R (E, B) = {\rm Hom}_R (E, R) \otimes_R B$$ for every
commutative $R$-algebra $B$.
\end{cor}

\begin{proof}
${\rm Hom}_R (E,B) = {\rm Hom}_R(E, R) \otimes_R B$ for every
commutative $R$-algebra if and only if ${\bf E}^*$ is a
quasi-coherent $R$-module. That is to say, if and only if ${\bf
E}$ is a module scheme, from the previous proposition, if and
only if $E$ is a projective module of finite type.
\end{proof}

\begin{cor}
An $R$-module $E$ is projective of finite type if and only if
$${\rm Hom}_R (E, E') = {\rm Hom}_R (E, R) \otimes_R E'$$ for every $R$-module $E'$.
\end{cor}

In \cite[Ch. II, \S 4.2, Prop. 2]{Bourbaki} it can be found the necessary
condition of this corollary. For the sufficient one, we will only
say that if ${\rm Id} = \sum_i w_i \otimes e_i \in {\rm
Hom}_R(E,E) = E^* \otimes_R E$, then $E = <e_i>$ is a module of finite type.

\begin{cor}
The quasi-coherent $R$-module corresponding to the $R$-module
$E^*$ is ${\bf E}^*$ if and only if $E$ is a projective module of
finite type.
\end{cor}


\begin{thebibliography}{99}


\bibitem[A]{Amel} \textsc{\'{A}lvarez, A., Sancho, C., Sancho,P.,}
\textit{\!Algebra schemes and their representations}, J.
Algebra {\bf 296/1} (2006) 110-144.

\bibitem[B]{Bourbaki} \textsc{Bourbaki, N.,}
\text{\!\'{E}l\'{e}ments de Math\'{e}matique. Alg\`{e}bre I}, Hermann, Paris, 1970.

\end{thebibliography}
\end{document}